\documentclass[a4paper,12pt,reqno]{amsart}
\usepackage{amsmath, amsthm, amscd, amsfonts, amssymb, graphicx, color}
\usepackage[bookmarksnumbered, colorlinks, plainpages]{hyperref}
\usepackage{wrapfig}
\usepackage{subfig}
 \usepackage[dvips]{epsfig}
\usepackage{geometry}
\usepackage{subfig}
\usepackage{cite}
\usepackage{esint}
\usepackage[framemethod=default]{mdframed}
\usepackage{showexpl}

\begin{document}
	
	 \newcommand{\be}{\begin{equation}}
	 \newcommand{\ee}{\end{equation}}
	 \newcommand{\bt}{\beta}
	 \newcommand{\al}{\alpha}
	 \newcommand{\laa}{\lambda_\alpha}
	 \newcommand{\lab}{\lambda_\beta}
	 \newcommand{\no}{|\Omega|}
	 \newcommand{\nd}{|D|}
	 \newcommand{\Om}{\Omega}
	 \newcommand{\h}{H^1_0(\Omega)}
	 \newcommand{\lt}{L^2(\Omega)}
	 \newcommand{\la}{\lambda}
	 \newcommand{\ro}{\varrho}
	 \newcommand{\cd}{\chi_{D}}
	 \newcommand{\cdc}{\chi_{D^c}}
      \def\avint{\mathop{\,\rlap{--}\!\!\int}\nolimits}

	 \newtheorem{thm}{Theorem}[section]
	 \newtheorem{cor}[thm]{Corollary}
	 \newtheorem{lem}[thm]{Lemma}
	 \newtheorem{prop}[thm]{Proposition}
     \newtheorem{assumption}{Assumption}
	 \theoremstyle{definition}
	 \newtheorem{defn}{Definition}[section]
	 \newtheorem{exam}{Example}[section]
	 \theoremstyle{remark}
	 \newtheorem{rem}{Remark}[section]
	 \numberwithin{equation}{section}
	 \renewcommand{\theequation}{\thesection.\arabic{equation}}
	 \numberwithin{equation}{section}
	 %
	 %
	 %
	\title[ Singularly
perturbed elliptic systems   ]{On a Class of Singularly Perturbed Elliptic Systems  with Asymptotic Phase Segregation}
	\author[Bozorgnia, Burger  ]{Farid Bozorgnia, Martin Burger}

 \address{Department of Mathematics, Instituto Superior T\'{e}cnico, Lisbon. } \email{bozorg@math.ist.utl.pt}

 \address{Department Mathematik, Friedrich-Alexander Universit\"at Erlangen-Nürnberg, Erlangen (FAU)   } \email{ martin.burger@fau.de }

	 \thanks{The corresponding author,   F. Bozorgnia was  supported by the Portuguese National Science Foundation through FCT fellowships SFRH/BPD/33962/2009}

	 \date{\today}

	
	 \begin{abstract}
	 	This work is devoted to  study of  a class of elliptic singular perturbed systems and their singular limit to a phase segregating system. We prove  existence and  uniqueness and study the asymptotic behaviour with convergence to a limiting problem   as  the interaction  rate  tends to infinity.  The limiting  problem is a free boundary problem such that at each point in the domain  at least one of the components is  zero which implies simultaneously all components can not coexist.   We present a novel method, which provides   an explicit solution of  limiting problem  for special   choice of parameters. Moreover, we present some numerical simulations of the asymptotic problem.

	 \end{abstract}

	 \maketitle

	 \textbf{Keywords}: Singular perturbed  system,  segregation,   free boundary problems, numerical approximation. \\
	 2010 MSC:58J37,35R35, 34K10.

\section{Introduction and problem setting}


In order to model   strong   interaction   between   multiple components  with     reaction and diffusion,  different models have been proposed. Among these models the adjacent segregation models have    been extensively studied  from different point  of views, to see about theoretical aspects  we refer to    \cite{CL, CR, DancerDu,  Ya}.   Most of the works  are   related to the case of two components, while \cite{CL} considers an extension to multiple components with strict segregation. Here we consider a different extension to multiple components that is still consist with the other models for the case of two components, the segregation behaviour is of different type for multiple ones however.

 Let $\Omega $ be bounded  domain with  $C^{1, \alpha}$   smooth boundary.   The model describes  the steady state   of $m$  species  diffusing  and   interacting  between all component   in   $\Omega.$  Let $u_{i}(x)$ denote the population density of the $i^{\textrm{th}}$  component.
We study  the following singular  elliptic system introduced in \cite{CR}, with unknowns  $U^{\varepsilon}=(u_{1}^{\varepsilon}, \cdots, u_{m}^{\varepsilon})$ which satisfy

\begin{equation}\label{s0}
\left \{
\begin{array}{llll}
\Delta  u_{i}^{\varepsilon}=   \frac{ A_{i}(x) }{\varepsilon}  F(u_{1}^{\varepsilon},\cdots, u_{m}^{\varepsilon}) & \textrm{ in  } \Omega,\\
u_{i}^{\varepsilon} \ge 0  & \textrm{ in  } \Omega,\\
u_{i} =\phi_{i} \,   &   \textrm{ on   } \partial \Omega,\\
 \end{array}
\right.
\end{equation}
 for  $i=1, \cdots, m$. Here   the function $F:\mathbb{R}^m \rightarrow \mathbb{R}$  is  given by
\[
F(u_1, \cdots ,u_m)=  \prod\limits_{j=1}^{m}  u_{j}^{\alpha_j},
\]
  for an  $m$-tuple $ (\alpha_1, \cdots, \alpha_m )$  with   $\alpha_i \ge 1.$

The  main assumptions on boundary values and data are as below:
\begin{assumption}

  The boundary data  $\phi_{i}$ are non-negative  $C^{1, \alpha}$  functions with  following partial segregation property
\[
\prod_{i=1}^{m} \phi_{i}=0 \quad   \textrm{on} \,\,  \partial \Omega.
\]
\end{assumption}

\begin{assumption}
  The functions $ A_{i}(x)$  are smooth, positive  and satisfy
 $$0<  A_{i}(x) \le \sum_{j\neq i}  A_{j}(x) \quad   \textrm{in} \, \,    \Omega.$$

\end{assumption}

  The  system (\ref{s0}) and  the limiting system for $\epsilon \downarrow 0$ appear  in   theory of flames   and are related to  a model     called  Burke-Schumann approximation.
 The main   assumption in  Burke-Schumann  model is that   oxidizer and reactant mix on a thin sheet and  the flame precisely occurs there.  A way to justify  the  underlying assumption is to
introduce a large parameter  called Damk\"{o}hler  number,   denoted by $D_a$, which is the parameter measuring the intensity of
the reaction  (see \cite{Will}). Then, the a chemical reaction is described by
\[
\textrm{Oxidizer}  +\textrm{ Fuel}   \rightarrow  \textrm{Products}.
\]
Let $Y_O$ and $Y_F$, respectively, denote the mass fraction of the oxidizer   and the fuel, then they  satisfy the following system
\begin{equation*}
\left \{
\begin{array}{llll}
-\Delta  Y_O + v(x).\nabla Y_O=   D_a \, Y_O \, Y_F  & \text{ in  } \Omega,\\
-\Delta  Y_F + v(x).\nabla Y_F=   D_a \, Y_O \, Y_F  & \text{ in  } \Omega,
 \end{array}
\right.
\end{equation*}
with given incompressible velocity field $v$ and a Dirichlet boundary condition on $\partial \Omega$.


In  \cite{CR}   a general H\"{o}lder estimate  for a   class of singular
perturbed elliptic system (\ref{s0})  is shown. The authors applied this estimate  to  the well-known Burke-Schumann
approximation in flame theory. Also they study the  classical cases  i,e., equidiffusional case with  high
activation energy approximation,  non- equidiffusional case, and to nonlinear diffusion
models.  The limiting problems are nonlinear elliptic equations; they have   H\"{o}lder   or Lipschitz
maximal global regularity.

We point out that L. Caffarelli   and   F. Lin  in \cite{CL}  studied the following   system with different coupling term
 \begin{align}\label{f20}
\begin{cases}
\Delta  u_{i}^{\varepsilon}=   \frac{ 1 }{\varepsilon}  u_{i}^{\varepsilon} \sum\limits_{j \neq i}   u_{j}^{\varepsilon} (x)\qquad\qquad & \text{ in  } \Omega,\\
u_{i}^{\varepsilon} \ge 0,\;  & \text{ in  } \Omega,\\
u^{\varepsilon}_{i}(x) =\phi_{i}(x)    &   \text{ on} \,  \partial \Omega,\\
  i=1,\cdots, m,
 \end{cases}
\end{align}
where the boundary values  satisfy
\[
\phi_{i}(x) \cdot \phi_{j}(x)=0,  \quad i \neq j \textrm{ on the boundary}.
\]

 \begin{rem}
 In system (\ref{s0})   choosing  $m=2$   and
 \[
 A_{i}(x)=1, \quad \alpha_{i}=1,   \,  i=1,2,
 \]
  we get   system (\ref{f20}) for $m=2$ which has been  studied extensively.  Thus  in  (\ref{s0}) we are interested when   $ m\ge 3. $

 \end{rem}

To see different theoretical aspects of the system  (\ref{f20})  we refer to  \cite{CL, Ya, W} and references therein. In   \cite{CL}   the authors  study the asymptotic limit;   as $\varepsilon $ tends to zero in system (\ref{f20}) and they show that limiting case     yields to pairwise segregation.  Furthermore, it is shown that  away from a  closed subset    of the Hausdorff dimension  less or equal $n-2 $  the free interfaces between various components   are, in fact, $C^{1,\alpha}$ smooth hyper surfaces.

For the  numerical approximation of the  system    (\ref{f20}) we refer to  \cite{BA,Bozorg1}. In  \cite{BA} the authors propose a numerical scheme  for a
class of  reaction-diffusion system with $m$  densities having disjoint supports and  are governed by a minimization problem. The proposed numerical scheme is applied
for the spatial segregation limit of diffusive Lotka-Volterra models in presence of high competition and inhomogeneous Dirichlet boundary conditions.  In \cite{AA} the proof of  convergence of the finite difference scheme for
a general class of  the spatial segregation  of reaction- diffusion, is given.

 This work is  devoted to analyse existence and uniqueness  results for system  (\ref{s0}), as well as a study of the qualitative properties of solutions to   (\ref{s0}) as $\varepsilon$ tends to zero. A particular novelty  of  the current     work is  to  provide an explicit solution for an arbitrary number of components $m$  when the parameter   $\varepsilon$ tends to  zero in the following system
  \begin{equation}\label{s1}
\left \{
\begin{array}{lll}
\Delta  u_{i}^{\varepsilon}=  \frac{A_{i}(x) }{\varepsilon}   \prod\limits_{j=1}^{m}  u_{i}^{\varepsilon}    & \text{ in  } \Omega,\\
u_{i}^{\varepsilon} \ge 0,  & \text{ in  } \Omega,\\
u_{i}(x) =\phi_{i}(x) \,     &   \text{ on   } \partial \Omega,
 \end{array}
\right.
\end{equation}
For the cases $A_{i}(x)$ be same or are constants.

 The outline of this paper is as follows:  Section 2    consists  the proof of existence and  uniqueness  of system (\ref{s0}). Section 3 deals with  the limiting case as  $\varepsilon$ tends to zero. In Section  4 we give an explicit solution for limiting case together with a rate of convergence.  Section 5 provides some numerical simulations of the singular limit.

\section{ Analysis  of the  model for fixed $\varepsilon$}

In this section we prove  existence and uniqueness  of the solution  of System  (\ref{s0})  for fixed   $\varepsilon$.   The proof is    constructive   and we implement it to obtain    numerical  approximation of  (\ref{s1})

 Consider the following related  time dependent parabolic   system
\begin{equation}\label{P1}
\left \{
\begin{array}{llll}
\frac{ \partial u_{i}^{\varepsilon}   }{\partial t} -   \Delta  u_{i}^{\varepsilon}=  - \frac{ A_{i}(x) }{\varepsilon}  F(u_{1}^{\varepsilon},\cdots, u_{m}^{\varepsilon}) & \text{ in  } \Omega \times (0, T)\\
u_{i}^{\varepsilon} (\cdot, 0)= u_{i0}   & \text{ in  } \Omega,\\
u_{i}(x,t) =\phi_{i}(x) \,   &   \text{ on   } \partial \Omega\times[0, T),\\
 \end{array}
\right.
\end{equation}
 where in (\ref{P1})  the initial values $ u_{i0}, \, i=1, \cdots m$ are non-negative and compatible with boundary data.  Then by Theorem 2.1 in \cite{EE} we obtain
 \[
 u_{i}^{\varepsilon}(x, t) \ge 0, \quad  t>0.
 \]
 Also it is straight to show that as $t$  tends to infinity
  \[
 u_{i}^{\varepsilon}(x, t) \rightarrow  u_{i}^{\varepsilon}(x),
 \]
 with $ u_{i}^{\varepsilon}(x)$ being the solution of (\ref{s0}), see \cite{CDH, MS}.

Let $(u_{1}^{\varepsilon}, \cdots  ,u_{m}^{\varepsilon})$  be a positive    solution of  the system  (\ref{s1}) then
\[
 u_{i} \le M \quad    i=1,\cdots ,m,
 \]
 where
\[
M= \underset{i=1, \cdots, m}{\max}\,\underset{ x\in \partial \Omega}{\max}\,  \phi_{i}(x).
\]
 We denote the harmonic extension of boundary data  $\phi_{i}$ with $ u_{i}^{0}$.  We multiply the following      equation
\[
\Delta (u_{i}^{\varepsilon}- u_{i}^{0}) =  \frac{A_i(x) }{\varepsilon}   \prod\limits_{j=1}^{m} F(u_{1}^{\varepsilon}, \cdots, u_{m}^{\varepsilon}).
\]
by $(u^{\varepsilon}_{i} - u_{i}^{0})^{+}$ where $u^{+}(x)=\max(u(x), 0).$
Then integrating by parts gives

\[
-\int_{\Omega} |\nabla( u_{i}^{\varepsilon}- u_{i}^{0}  )^{+} |^{2} dx= \int_{\Omega}\frac{A_i(x) }{\varepsilon} ( u_{i}^{\varepsilon} - u_{i}^{0})^{+}  F(u_{1}^{\varepsilon}, \cdots, u_{m}^{\varepsilon})\, dx.
\]
 Note that the integrand of right hand side is positive and
 \[
 u_{i}^{\varepsilon}-  u^{0}_{i}= 0 \quad  \text{on }  \partial\Omega.
 \]
 From here  $ \int_{\Omega} |\nabla( u_{i}^{\varepsilon}- u_{i}^{0}  )^{+} |^{2} dx= 0$ which implies
 \[
 u_{i}^{\varepsilon} \le  u_{i}^
{0} \quad \textrm{in} \, \Omega.
 \]

  A standard maximum and nonnegativity principle for elliptic equations (cf. \cite{schaefer}) yields the following result. In sequel we use this result.

\begin{lem}\label{sysn1}
Let $u  \in H^1(\Omega)$  be a  weak  solution of  the system
 \begin{equation}
\left \{
\begin{array}{lll}
\Delta  u =    a  \, u^{\alpha}(x)   & \text{ in  } \Omega,\\
u  =\phi   \,     &   \text{ on   } \partial \Omega.
 \end{array}
\right.
\end{equation}
with $a $ and $\phi$ bounded and nonnegative, $\alpha\geq 1$
then
\[
0\le u \le M,
\]
where
\[
M=  \underset{ x\in \partial \Omega}{\max}\,  \phi_{i}(x).
\]
\end{lem}

In  the next   Theorem \ref{sun0} we show   the  existence of nonnegative    solutions to the original system.  The main idea of the proof is to construct sub and super solution  and decoupling the system in iterative  way and to exploit the uniform $L^\infty$ bounds, see also  the proof in \cite{W} for  the proof of uniqueness of the solution for  system (\ref{f20}).

  \begin{thm}\label{sun0}
For each $\varepsilon >0, $ there exist a unique nonnegative solution
$$(u_{1}^{\varepsilon},\cdots, u_{m}^{\varepsilon}) \in H^1(\Omega)^m \cap L^\infty(\Omega)^m$$ of the system  (\ref{s0}).
 \end{thm}

 \begin{proof}
 Without loss of generality  in the proof we set $  \alpha_{i} =1$ i.e.,
  \[
F(u_1, \cdots ,u_m)=  \prod\limits_{j=1}^{m}  u_{j}.
\]
 To start, consider the harmonic extension $u_{i}^{0}$ given by
 \begin{equation}\label{sys2}
 \left \{
 \begin{array}{llll}
  - \Delta u_{i}^{0} = 0   & \text{ in  } \Omega,\\
   u_{i}^{0}  =\phi_{i}     &   \text{on   } \partial \Omega.
  \end{array}
 \right.
 \end{equation}
Next,  given $ u_{i}^{k}$ consider the solution of the following linear   system
 \begin{equation}\label{sy4}
 \left \{
 \begin{array}{lllll}
 \Delta  u_{i}^{k+1}=\frac{A_{i}(x)}{\varepsilon}  \, \frac{
 u_{1}^{k} \cdots u_{i-1}^{k}u_{i}^{k+1} u_{i+1}^{k} \cdots u_{m}^{k} \, + \,  u_{1}^{k+1} \cdots u_{i-1}^{k+1} u_{i}^{k+1} u_{i+1}^{k} \cdots u_{m}^{k}}{2}        & \text{ in  } \Omega,\\
 u_{i}^{k+1}(x) =\phi_{i}(x)    &   \text{ on   } \partial \Omega.\\
  \end{array}
 \right.
 \end{equation}
Note that we can subsequently solve the equations for increasing $i$ due to the triangular structure and always obtain a problem of the form considered in Lemma \ref{sysn1}, hence the uniform bounds apply.
 We show that the following inequalities hold:
 \[
 u_{i}^{0}\ge u_{i}^{2} \cdots \ge u_{i}^{2k}\ge \dots  \ge u_{i}^{2k+1}\ge \cdots \ge u^{3}_{i}\ge u_{i}^{1},  \quad \textrm{in} \, \Omega.
 \]
 The first iteration for $u_{1}$ reads as
 \[
 \Delta  u_{1}^{1}=  \frac{A_{1}(x)}{\varepsilon}  u_{1}^{1}  u_{2}^{0}\cdots u_{m}^{0}.
 \]

Note that since $u_{i}^{0} \ge 0, $ and    boundary conditions  $\phi_{i}(x)$  are non negative then the weak maximum principle (see appendix)  implies  that $ u_{1}^{1}\ge 0.$
The  equation for $u_{2}^{1}$ in (\ref{sy4}) is given by
\[
 \Delta  u_{2}^{1}=  \frac{A_{2}(x)}{2\varepsilon}\,  (\, u_{1}^{0} u_{2}^{1}\,   u_{3}^{0}\cdots u_{m}^{0}  \, +  \,  u_{1}^{1} u_{2}^{1}\,   u_{3}^{0}\cdots u_{m}^{0} \,  ).
 \]
 Repeating the same argument,  we obtain that  $u_{2}^{1}\ge 0$   and     consequently
\[
 u_{i}^{1}\ge 0, \quad i=3, \cdots ,m.
 \]
 Now we have
  \begin{equation}\label{sy44}
  \left \{
  \begin{array}{ll}
  \Delta  u_{i}^{1}\ge 0 & \text{ in  } \Omega,\\
  u_{i}^{1}(x) =u_{i}^{0}(x)= \phi_{i}(x)    &   \text{ on   } \partial \Omega.\\
   \end{array}
  \right.
  \end{equation}
Thus  the comparison principle  implies that $ u_{i}^{1}\le  u_{i}^{0} $.  The same argument shows
 \[
 u_{i}^{0} \ge u_{i}^{2}.
 \]
 In the next step we verify the following inequalities hold
 \[
u_{i}^{2} \ge u_{i}^{1}\quad i=1, \cdots,m.
 \]
 To do this, one verifies that  inequality   $u_{1}^{2} \ge u_{1}^{1}$ holds  then   this fact  can be used to prove inequality  for  $i=2,3,\cdots,m$.  Then the same arguments  show that
 \[
u_{i}^{3} \ge u_{i}^{1}.
 \]
 To proceed more with induction, assume that
  \begin{equation}\label{inq2}
 u_{i}^{0}\ge u_{i}^{2}\ge  \cdots \ge u_{i}^{2k}\ge  u_{i}^{2k+1}\ge \cdots \ge u^{3}_{i}\ge u_{i}^{1}.
  \end{equation}
 We show that
 \[
 u_{i}^{2k+1}\le u_{i}^{2k+2}.
 \]
 To show this,  first we  check  for $i=1$ and the same argument can be applied consequently.
  By (\ref{sy4}) and the  assumption in  (\ref{inq2})  we have

   \begin{equation*}
  \left \{
  \begin{array}{ll}
   \Delta  u_{1}^{2k+2}= \frac{A_{1}(x) }{\varepsilon}  u_{1}^{2k+2} \prod\limits_{j=2}^{m}u_{j}^{2k+1} \le  \frac{ 1 }{\varepsilon}  u_{1}^{2k+2} \prod\limits_{j=2}^{m}u_{j}^{2k},\\\\
    \Delta  u_{1}^{2k+1}=  \frac{A_{1}(x)}{\varepsilon}  u_{1}^{2k+1} \prod_{j=2}^{m}u_{j}^{2k}.
   \end{array}
  \right.
  \end{equation*}

 Note that  $  u_{1}^{2k+1}  $ and $  u_{1}^{2k+2}$ have the same boundary value so by the comparison principle
 \[
  u_{1}^{2k+1}\le u_{1}^{2k+2}.
  \]
 Now we proceed for $i=2,\cdots,m$.
 The same argument  using the assumption  $u_{i}^{2k+1} \ge u_{i}^{2k-1}$   shows that
   \[
    u_{i}^{2k+2}\le u_{i}^{2k}.
    \]
   For the next step,  we use the fact from previous step  which states  $  u_{i}^{2k+2}\le u_{i}^{2k}$  to verify
   $  u_{i}^{2k+3}  \ge  u_{i}^{2k+1}. $

   Now let $ \overline{u}_{i} $ and  $ \underline{u}_{i} $ be two families of functions such that
   \[
    u_{i}^{2k} \rightarrow   \overline{u}_{i} \quad  \textrm{uniformly in } \Omega,
    \]
     \[
        u_{i}^{2k+1} \rightarrow \underline{u}_{i} \quad  \textrm{uniformly in } \Omega.
        \]
Taking the limit in   (\ref{sy4}) yields  for $i=1,\cdots m$ the followings hold

  \begin{equation}\label{s9}
  \left \{
  \begin{array}{llll}
  \Delta \overline{u}_{i} =  \frac{A_{i}(x)}{2\varepsilon}( \overline{u}_{1}\cdots    \overline{u}_{i} \underline{u}_{i+1}\cdots \underline{u}_{m}  +\underline{u}_{1}\cdots \underline{u}_{i-1} \overline{u}_{i} \underline{u}_{i+1} \cdots \underline{u}_{m})     & \text{ in  } \Omega,\\
 \Delta \underline{u}_{i} =  \frac{ A_{i}(x) }{2\varepsilon}( \underline{u}_{1}\cdots   \underline{u}_{i}   \overline{u}_{i+1}\cdots   \overline{u}_{m} +   \overline{u}_{1}\cdots \overline{u}_{i-1}\underline{u}_{i} \overline{u}_{i+1}\cdots  \overline{u}_{m} )    & \text{ in  } \Omega.
       \end{array}
  \right.
  \end{equation}
 The inequality $    u_{i}^{2k}   \ge u_{i}^{2k+1}$ implies that
 \begin{equation}\label{in1}
  \overline{u}_{i} \ge \underline{u}_{i}.
 \end{equation}
   We will show that in fact the equality holds.  To  do this, first   consider the equations for the $m^{\textrm{th}}$
   \begin{equation}\label{s10}
  \left \{
  \begin{array}{llll}
  \Delta \overline{u}_{m} =  \frac{ A_{m}(x) }{2\varepsilon}\, \overline{u}_{m}\left( \overline{u}_{1}\cdots    \overline{u}_{i} \, \overline{u}_{i+1} \cdots \overline{u}_{m-1} +  \underline{u}_{1}\cdots   \underline{u}_{i}  \underline{u}_{i+1}  \cdots   \underline{u}_{m-1} \right)    & \text{ in  } \Omega,\\
 \Delta \underline{u}_{m} =  \frac{A_{m}(x) }{2\varepsilon}\,\underline{u}_{m}  \left(\underline{u}_{1}\cdots   \underline{u}_{i}  \underline{u}_{i+1}  \cdots   \underline{u}_{m-1} \,   +  \, \overline{u}_{1}\cdots    \overline{u}_{i} \, \overline{u}_{i+1} \cdots \overline{u}_{m-1} \right)   & \text{ in  } \Omega,\\
\overline{u}_{m}= \underline{u}_{m}=  \phi_{m}(x)    &   \text{ on   } \partial \Omega,
       \end{array}
  \right.
  \end{equation}
   which implies
    \[
       \overline{u}_{m}  = \underline{u}_{m}.
       \]
           Now by checking the equation for $i=m-1$  in (\ref{s9})  and using the previous  fact
               $\overline{u}_{m}  = \underline{u}_{m},$  yields   \[
  \overline{u}_{m-1}  = \underline{u}_{m-1},
  \]
   and   argument is repeated backward  which  shows equality for every $i$.

 To show uniqueness,    assume there exists  another positive  solution $(w_1,\cdots,w_m)$ of  system, then we show
\[
u_{i}=w_{i}, \quad i=1, \cdots ,m.
\]
We will prove that the following equations hold:
 \begin{equation}\label{ineq1}
u_{i}^{2m+1} \le w_{i}\le u_{i}^{2m} , \quad  \textrm{ for } \,    m \ge 0.
 \end{equation}
To begin, we show that
 \begin{equation}\label{ineq20}
w_{i}\le u_{i}^{0}.
\end{equation}
This is a consequence of the fact that   $w_{i}$ satisfies
 \begin{equation*}
  \left \{
  \begin{array}{llll}
   \Delta w_{i}\ge 0    & \text{ in  } \Omega,\\
     w_{i}=u^{0}_{i}   & \text{on   } \partial \Omega.
       \end{array}
  \right.
  \end{equation*}
Next we compare $w_{i}$  with $u_{i}^1$ and we show $ w_{i} \ge u_{i}^1$. As in existence part, first we check for $i=1$ in inequality follows from (\ref{ineq20}) and
 \begin{equation*}
  \left \{
  \begin{array}{ll}
   \Delta w_{1} =\frac{ A_1 w_1}{\varepsilon}  \prod\limits_{j=2}^{m}  w_j   & \text{ in  } \Omega,\\
    \Delta u_{1}^1   =\frac{A_1 u_{1}^{1}}{\varepsilon}   \prod\limits_{j=2}^{m}  u_{j}^{0}  & \text{ in  } \Omega.
       \end{array}
  \right.
  \end{equation*}
 Now we proceed by induction and we assume that the claim is true until $2k+1.$  This means that we have
\[
u_{i}^{2k+1}\le w_{i} \le u_{i}^{2k}.
\]
Then we show
\[
u_{i}^{2k+3}\le w_{i} \le u_{i}^{2k+2}.
\]
Again comparing  the equations  for $w_i$ and $u_{i}^{2k+2} $  and the using assumption $ u_{i}^{2k+1} \le w_i$ yields
the following inequality
\[
w_{i} \le u_{i}^{2k+2}.
\]%
 The same reasoning for  inequality $ u_{i}^{2m+3}\le   w_{i}$  holds.  Now taking limit in (\ref{ineq1} ) shows that
  \[
  w_i=u_i, \quad i=1, \cdots,m.
  \]
\end{proof}


\section{Limiting problem}

In this section we  study properties of the solution for  system  (\ref{s0})  to    provide estimates  and compactness results  to pass to the limit  as $\varepsilon $ tends to zero.

As  we have seen in  the last section,  for each fixed  $\varepsilon,$ the system (\ref{s0}) has a  unique solution.   Let $ U^{\varepsilon}=(u_{1}^{\varepsilon},  \cdots ,u_{m}^{\varepsilon} )$ be the unique  positive  solution of system (\ref{s0}) for fixed $\varepsilon,$ then   $u_{i}^{\varepsilon} $  for $i=1, \cdots,m$  satisfy  the following differential inequalities:
\[
 -\Delta u_{i}^{\varepsilon}  \le 0 \quad \text{ in } \quad \Omega.
 \]
 Also  define $ \widehat{u}_{i}^{\varepsilon}$  as
  \begin{equation}\label{hat}
  \widehat{u}_{i}^{\varepsilon}  := u_{i}^{\varepsilon}  - \sum_{j\neq i}u_{j}^{\varepsilon},
   \end{equation}
   then considering the assumption 2, it is  easy to verify that
 $$  -  \Delta    \widehat{u}_{i}^{\varepsilon}  \ge 0. $$
 Let   $h_{i}$ and   $  H _{i}$  for $i=1, \cdots m $  be harmonic with boundary value  $\phi_{i} $  and
 $\widehat{\phi}_{i}$ respectively, where
\[
  \widehat{\phi}_{i}=\phi_{i}- \sum_{j\neq i} \phi_{j},
  \]
  then we have
  \[
H_{i} \le   \widehat{u_i}^{ \varepsilon} \le     u_{i}^{\varepsilon}  \le h_{i},
\]
which implies
 \begin{equation}\label{EE}
\frac{\partial h_{i}}{\partial \nu} \le
    \frac{\partial u_{i}^{\varepsilon}}{\partial \nu}.
\end{equation}

 In this part we show that the solution  $ u_{i}^{\varepsilon}$  of  system  (\ref{s0}) has bound in $W^{1,2}(\Omega) $ independently of $\varepsilon.$  To do this, we prove several lemmas.
\begin{lem}\label{sunf1}
Assume $x_0\in \Omega $ and $ B_{2r}(x_0) \subset \Omega$.  Let $u$ satisfies the following
 \begin{equation*}
\left \{
\begin{array}{lll}
\Delta u= f \ge 0   & \text{ in  }\,  B_{2r}(x_0),\\
 0\le u \le M    &   \text{in   } \,  B_{2r}(x_0).
 \end{array}
\right.
\end{equation*}
  Then
\[
\int_{B_{r}(x_0)}  f(x) \, dx \le C_{0} M r^{n-2},
\]
 for some $C_{0}$ that only depends on  dimension $n$.
 \end{lem}
 \begin{proof}
   Without loss of generality, assume $x_0=0$. By Green's formula for ball one has
   \[
   0\le u(0)=\avint_{\partial B_{2r}(0)}   u(x) \, dx -  \int_{B_{2r}(0)} ( \frac{\omega_n}{|x|^{n-2}}- \frac{\omega_n}{(2r)^{n-2}}) f(x) \, dx
   \]
   \[
   \le M - C_{0} \int_{B_{r}}\frac{  \omega_n}{|x|^{n-2}}f(x)\, dx \le M -\frac{C_{0} }{r^{n-2}}\, \omega_n \int_{B_{r}(0)}f(x)\, dx.
   \]
  Next, rearranging terms proves the Lemma.
 \end{proof}
\begin{lem}\label{sunf2}
Assume that $u_{i}$ satisfies
\begin{equation}\label{sunf3}
\left \{
\begin{array}{lll}
\Delta u_{i}= \frac{1}{\varepsilon}   \prod\limits_{j=1 }^{m}  \,  u_{j}  & \text{ in }\ \Omega, \\
      u_{i} \ge 0  & \text{ in } \,  \Omega,\\
      u_{i}  =\phi_{i}   &  \text{ on } \  \partial\Omega.
\end{array}
\right.
\end{equation}
Then there exists a constant  $C_0$ depends only on $\Omega, n, r$ and
 $\|\phi_{i}\|_{C^{1, \alpha}(\partial \Omega)} $ such that
\[
\fint_{B_{r}(x_0)  \cap  \Omega} \frac{1}{\varepsilon}   \prod\limits_{j=1 }^{m}  \,  u_{j}\, dx \le C_{0} r^{n-2}.
\]
\end{lem}
\begin{proof}
The proof consider different cases.
\begin{enumerate}
  \item  If  $B_{2r}(x_0) \subset \Omega$ then it follows  by previous  Lemma \ref{sunf1}.
  \item  If   $ \exists k $ such that $ \phi_{k}=0 $ on $ \partial\Omega \cap  B_{2r}(x_0)$ then we may extend  $u_k$ to
      \begin{equation*}
\overline{u}_{k}= \left \{
\begin{array}{ll}
 u_k  & \text{ in }\ \Omega, \\
   0   &  \text{ in  } \,   \Omega^{c}.
    \end{array}
\right.
\end{equation*}
      and apply the previous Lemma to $\overline{u}_k.$
  \item If none of $\phi_{k}$  vanishes on $ \partial\Omega \cap  B_{2r}(x_0)$ then, since the product of boundary values is zero,  there must be a $ \phi_{i}$ that vanishes at a point
      $ y_1\in  \partial\Omega \cap  B_{2r}(x_0),$  we may assume that  $ \phi_{1}(y_1)=0.$   Also, since
      $ u_1 \ge  0 $ it follows that
      \[
\frac{\partial u_{1}(y_1)}{\partial \nu} \le 0.
      \]
Now let $h_1$ solves
\begin{equation}\label{sunf4}
\left \{
\begin{array}{lll}
\Delta h_1=  0 & \text{ in }\,  B_{4r}(x_0) \cap  \Omega, \\
  h_1  =\phi_{1}   &  \text{ on } \  \partial\Omega  \cap  B_{4r}(x_0),\\
   h_1= 0  & \text{ on} \,  \Omega \cap \partial   B_{4r}(x_0).
\end{array}
\right.
\end{equation}
Since  $\partial \Omega$ and $ \phi_{1}$ are $C^{1, \alpha}$ it follows that
\[
|\nabla h_1 | \le C_h, \quad  \text{ in  } \,  \Omega  \cap  B_{3r}(x_0).
\]
Now either $w=u_1- h_1$ satisfies
\[
-C_h \le \frac{\partial w  }{\partial \nu} \le C_{*},
\]
 for some $C_{*}$ (which we will decide ) in which case  we may apply the previous Lemma on $w$ since
\[
\Delta w= \frac{1}{\varepsilon}   \prod\limits_{j }^{m}  \,  u_{j} + \frac{\partial w }{\partial \nu} \mathcal{H}^{n-1}\lfloor_{\partial\Omega} \quad \textrm{  in} \,  B_{3r}(x_0),
\]
or there is a point $y_2 \in  \partial \Omega  \cap  B_{3r}(x_0)$ such that
 \[
\frac{\partial w(y_2)}{\partial \nu} \ge C_{*}.
      \]
      Note that $ \phi_{1} (y_2)  > 0$ since otherwise,
      \[
      0 \ge \frac{\partial u_{1}(y_2)}{\partial \nu}=  \frac{\partial w(y_2)}{\partial \nu}+ \frac{\partial h_{1}(y_2)}{\partial \nu} \ge C_{*} - C_h > 0,
      \]
      provided  $C_*$ is large enough.  Next, since $ \phi_{1} (y_2)> 0$ there is another $ \phi_{k}$ say
      $ \phi_{2}$, such that $ \phi_{2} (y_2)  =0.$   Let $h_2$ solves

      \begin{equation}\label{sunf5}
\left \{
\begin{array}{lll}
\Delta h_2=  0 & \text{ in }\,  B_{4r}(x_0) \cap  \Omega \\
  h_2  =\phi_{2}   &  \text{ on } \  \partial\Omega  \cap  B_{4r}(x_0),\\
   h_2= 0  & \text{ on} \,  \Omega \cap \partial   B_{4r}(x_0).
\end{array}
\right.
\end{equation}
 Then again $ |\nabla h_2 | \le C_h$ in $  B_{3r}(x_0)$ for some $C_h$ depending only on the domain $\Omega$ and   $\|\phi_{2}\|_{C^{1, \alpha }}$.
 Next let  $ u_2= w + h_2+ g$ in  $B_{4r}(x_0)\cap \Omega  $  where

   \begin{equation}\label{sunf6}
\left \{
\begin{array}{lll}
\Delta g=  0 & \text{ in }\,  B_{4r}(x_0) \cap  \Omega \\
 g =u_2 -w   &  \text{ on } \  \partial\Omega  \cap  B_{4r}(x_0),\\
   g= 0  & \text{ on} \,    \partial   B_{4r}(x_0)  \cap\Omega.
\end{array}
\right.
\end{equation}

Since $g$ is bounded;  $|g|\le 3 M$    on $\partial   B_{4r}(x_0)  \cap\Omega,$  then    it follows that
\[
 |\nabla g| \le C_g, \quad \text{  in  }   B_{3r}(x_0) \cap \Omega,
\]
where $C_g$ depends on the bound $ M, r $ and $\Omega$. This leads in particular to
 \[
      0 \ge \frac{\partial u_{2}(y_2)}{\partial \nu}=  \frac{\partial w(y_2)}{\partial \nu}+ \frac{\partial h_{1}(y_2)}{\partial \nu} +  \frac{\partial g(y_2)}{\partial \nu}  \ge C_{*} - C_h -C_g \ge  0.
      \]
This is a contradiction if $C_{*}$ is large enough and this complete the proof.
\end{enumerate}
\end{proof}
\begin{prop}\label{sunf16}
Let $ u_{1} , \cdots ,u_m$ be as in previous Lemma. Then there exists a constant $C_0$(independent of $\varepsilon$ such that
\[
\| u_i\|_{W^{1,2}(\Omega)} \le  C_0.
\]
\end{prop}
\begin{proof}
  Cover $\Omega$ by finitely say $N$  balls $B_{r}(x_k)$ and notice that
  \[
  \int_{\Omega} \frac{1}{\varepsilon}   \prod\limits_{j }^{m}  \,  u_{j} \le  \sum_{k=1}^{N} \int_{B_{r}(x_k) } \frac{1}{\varepsilon}   \prod\limits_{j }^{m}  \,  u_{j} \le N\, C_{0} r^{n-2}.
  \]
  Next let $f=\frac{1}{\varepsilon}   \prod\limits_{j=1 }^{m}  \,  u_{j} $  and define
  \[
  v= -\int_{\Omega} \frac{\omega_n}{|x-y|^{n-2} }f(y) \, dy.
  \]
  Then $v$ satisfies
  \[
  \Delta v= f \chi_{\Omega} \quad \text{in} \,\mathbb{R}^n,
    \]
    which implies
     \begin{equation}\label{sunf7}
\int_{B_{R}(0)} |\nabla v|^{2}\,  dx = \int_{B_{R}(0)} f(x) v(x) \,  dx+ \int_{\partial B_{R}(0)}      v\, \frac{\partial v}{\partial \nu}\, ds  \le C.
 \end{equation}
    where $R$ is chosen so large that $\Omega  \subset B_{R}(0).$   Now let $ u_i= H_i +v$ where
    \begin{equation}\label{sunf8}
\left \{
\begin{array}{lll}
\Delta H_i=  0 & \text{ in }\,    \Omega \\
  H_i  =\phi_i- v   &  \text{ on } \  \partial\Omega.
\end{array}
\right.
\end{equation}
Since  $ \phi_i \in C^{1, \alpha} $ and $ v\in W^{1,2}(\Omega)$ by  (\ref{sunf7}), then it follows that
$ H_i  \in W^{1,2}(\Omega) $ with bounds only depending on $ \| v\|_{W^{1,2}(\Omega)},   \| \phi_i\|_{C^{1,\alpha}(\Omega)}  $  and $(\Omega). $  In particular, $ \| u_i\|_{W^{1,2}(\Omega)}$ is bounded independent of $\varepsilon $.
\end{proof}

The above Lemma  shows that   up to a subsequence  denoted with $u_{i}^{\varepsilon} $ we get
   \[
   u_{i}^{\varepsilon} \rightharpoonup  u_{i} \quad \text{in} \quad H^{1}_{0}(\Omega).
   \]

 The main result of this section is  Theorem \ref{lim}  which shows the asymptotic behaviour of system (\ref{s1})
 as  $\varepsilon$ tends to zero.

 \begin{thm}\label{lim}
     Let $U^{\varepsilon}=(u_{1}^{\varepsilon},\cdots ,u_{m}^{\varepsilon})$ be a solution of the  system    at fixed $\varepsilon$.  Let $ \varepsilon $ tends to zero, then there exists $ U \in  H^{1}(\Omega))^{m}  \cap L^\infty(\Omega)^m$ such
that for all $ i=1,\cdots,m$:
\begin{enumerate}

\item $\Delta u_i \geq 0$  in the sense of distribution.

 \item up to  subsequences, $u_{i}^{\varepsilon}- u_{i} \rightarrow  0$  strongly in $H^{1}_{0} (\Omega)$.

  \item $\prod\limits_{i}^{m} u_{i}=0$    a.e in \ $ \Omega.$

\end{enumerate}
\end{thm}

\begin{proof}

Proposition  (\ref{sunf16})  shows  the existence of a weak limit $u_{i}$  such that, up to subsequences,
\[
u_{i}^{\varepsilon}\rightharpoonup u_{i} \quad  \textrm{in  } \, H_{0}^{1}.
\]
The weak limit $u_{i}$  for $i=1, \cdots, m$ satisfy the following differential inequalities
 \[
  -\Delta u_{i} \le 0,  \quad    -  \Delta    \widehat{u}_{i}    \ge 0        \text{ in } \quad \Omega,
 \]
 since  we can pass to the weak  limit in  the differential inequalities  for  $ u_{i}^{\varepsilon}$ and  $ \widehat{u}_{i}^{\varepsilon}.$   To show the strong convergence, we  show that
  \[
  \int_{\Omega}  |  \nabla u_{i}^{\varepsilon} |^{2} \,   dx \rightarrow  \int_{\Omega}  | \nabla u_{i}|^{2}\,   dx.
  \]
   By weak lower semi continuously  of Dirichlet norm  just needs to show
   \[
   \int_{\Omega}  | \nabla  u_{i}|^{2} \,   dx  \ge   {\lim \sup}  \int_{\Omega}  | \nabla u_{i}^{\varepsilon}  |^{2} \,   dx.
   \]
We multiply the inequality  $-\Delta  u_{i}^{\varepsilon} \le  0 $  by  $u_{i}^{\varepsilon}$   and integration by parts,
  \[
   \int_{\Omega}    |\nabla u_{i}^{\varepsilon}|^{2} \,   dx   -    \int_{\partial \Omega}   u_{i}^{\varepsilon}  \frac{\partial u_{i}^{\varepsilon}}{\partial n}  \,   ds   \le 0.
   \]
   This implies
    \begin{equation}\label{1E}
           \int_{\partial \Omega}   u_{i}  \frac{\partial u_{i} }{\partial n}  \,   ds  \ge  {\lim \sup}  \int_{\Omega}    |\nabla u_{i}^{\varepsilon}|^{2} \,   dx.
    \end{equation}
Next we multiply the equation for  $u_{i}^{\varepsilon} $ by $ u_{i}$ to obtain
\[
 - \int_{\Omega}     \nabla u_{i}^{\varepsilon} \cdot  \nabla u_{i}   \,   dx  +  \int_{\partial \Omega}   u_{i}   \frac{\partial u_{i}^{\varepsilon}}{\partial n}  \,   ds  = \int_{\Omega}  u_{i} \,  \prod\limits_{j=1}^{m}  u_{j}^{\varepsilon} \, dx.
   \]
Taking the limit as  $\varepsilon_{n}$ tends to zero and considering the weak convergence of $u_{i}^{\varepsilon}$ and previous part to have
 \begin{equation}\label{2E}
 - \int_{\Omega}      |  \nabla u_{i} |^{2}  \,   dx  +  \int_{\partial \Omega}   u_{i}   \frac{\partial u_{i}}{\partial n}  \,   ds  = 0.
    \end{equation}
Form (\ref{1E}) and (\ref{2E}) the result holds.

 \item  (2)
Fix a point $x_{0} \in \Omega$ and  let the index $i$ be such that
  \[
  u^{\varepsilon}_{i}(x_0)=\underset{ 1\le k\le m}  {\max} u^{\varepsilon}_{k}(x_{0}).
  \]
  Now assume that $u_{i}^{\varepsilon}(x_{0})=c > 0 $ then by H\"{o}lder continuity there is $r$ such that
 \[
 |u_{i}^{\varepsilon}(x)- u_{i}^{\varepsilon}(x_{0})| \le \frac{c}{2}, \quad x \in B(x_{0}, r).
 \]
 Next we use the fact that the functions  $u_{i}^{\varepsilon}$  for $i=1, \cdots,m$ are subharmonic, using the mean value property for subharmonic functions (see the proof of theorem 2.1 in \cite{GT})

\begin{equation} \label{eq_sun}
 \begin{split}
 \avint_{\partial B(x_{0},r)}  | u_{i}^{\varepsilon}(x_{0}) - u_{i}^{\varepsilon}(y) |  \, dy   & = \int_{0}^{r} ( \int_{ B(x_{0},s)}  \Delta u_{i}^{\varepsilon} ) \, \frac{ds}{s^{n-1}}  \\
 & \ge r^{2} \avint_{B(x_{0},r)}  \Delta u_{i}^{\varepsilon} \, dx.
 \end{split}
\end{equation}
 From here the following holds
 \begin{equation}\label{E1}
 \avint_{B(x_{0},r)}  \Delta u_{i}^{\varepsilon} \, dx =  \avint_{B(x_{0},r)}\frac{u_{i}^{\varepsilon}}{\varepsilon}    \prod\limits_{j\neq i}^{m}  u_{j}^{\varepsilon}  (x)\le \frac{c}{2r^2}.
 \end{equation}
 Note that in the  ball $B(x_{0},r)$ we have $u_{i}^{\varepsilon}  \ge \frac{c}{2}   $  so from (\ref{s1})  we obtain
\begin{equation}\label{sunshine2}
 \avint_{B(x_{0},r)}\frac{1}{\varepsilon}    \prod\limits_{j\neq i}^{m}  u_{j}^{\varepsilon}  (x)\, dx \le \frac{1}{r^2}.
\end{equation}
Next, in (\ref{sunshine2})   let  $\varepsilon$  tend to zero  which yields
\[
\prod\limits_{j\neq i}^{m}  u_{j}^{\varepsilon}  (x) \rightarrow 0 \quad \text{in} \quad B(x_{0},r).
\]

\end{proof}

 Let $w_1$ be the  first eigenfunction of the Laplace operator in $\Omega,$   i.e.,
 \begin{equation*}
  \left \{
  \begin{array}{ll}
   -\Delta w_{1}=\lambda_{1} w_1    & \text{ in  } \Omega,\\
     w_1=0    & \text{on   } \partial \Omega.
       \end{array}
  \right.
  \end{equation*}
  The first eigenfunction does not change the sign and  we may therefore take it to be positive and  normalized it so that
    $\| w_1\|_{L^\infty} =1.$    Multiplying the equation
\[
\Delta  u_{i}^{\varepsilon}=   \frac{ A_i(x) }{\varepsilon}   \prod\limits_{i=1}^{m}  u_{i}^{\varepsilon}  (x),
\]
by $w_1$ and integrating over $\Omega$ yields
\[
\int_{\Omega} w_1\,   \Delta  u_{i}^{\varepsilon}\,   dx = \int_{\Omega}  \frac{ A_i(x) }{\varepsilon}   \prod\limits_{i=1}^{m}  u_{i}^{\varepsilon}  (x) \, w_1 \, dx.
\]
Integration by parts and implementing that $w_1$ is zero  on boundary, we obtain
\begin{equation*}
\begin{split}
\int_{\Omega}  \frac{ A_i(x) }{\varepsilon}   \prod\limits_{i=1}^{m}  u_{i}^{\varepsilon}  (x) \, w_1 \, dx & = \int_{\Omega} u_{i}^{\varepsilon}\, \Delta w_1  \,  dx-   \int_{ \partial \Omega} u_i \,  \frac{ \partial w_1 }{\partial n}\, ds  \\
 & = \lambda_{1} \int_{\Omega} u_{i}^{\varepsilon}\,  w_1  \,  dx-   \int_{ \partial \Omega} \phi_{i}  \,  \frac{ \partial w_1 }{\partial n}\, ds.
\end{split}
\end{equation*}
 Now  from the bound on  $u_i$ and  the fact that normal derivative of the first eigenfunction  on the boundary is bounded, we conclude
 \[
 \int_{\Omega}  \frac{ A_i(x) }{\varepsilon}   \prod\limits_{i=1}^{m}  u_{i}^{\varepsilon}  (x) \, w_1 \, dx \le C.
 \]
  We know that  for $i=1, \cdots ,m$ the solution $u_{i}^{\varepsilon}$  are H\"{o}lder continuous
 \[
\| u_{i}^{\varepsilon}\|_{C^{\alpha}} \le C_i,
 \]
   where  the constant  $C_i$ is independent   of $\varepsilon.$
  Note that,  since
 \[
 w_{1}(x)  > 0 \quad  \textrm{in the interior of }  \Omega,
  \]
  then the inequality in above yields that
  \begin{equation}\label{bound1}
 \int_{\Omega'}  \Delta  u_{i}^{\varepsilon} \le C \quad   \text{ in  compact subsets  }\Omega' \subset \Omega,
  \end{equation}
  where the constant $C$ is independent of  $ \varepsilon.$  For the rest, we show that for those  points close  enough
   to the boundary, $\Delta  u_{i}^{\varepsilon}$   remains bounded.   Let $C_{i}$   and   $\beta_{i} $ denote  the H\"{o}lder constant and   H\"{o}lder exponent  of $u^{\varepsilon}_{i}.$    Choose the strip around boundary such that
\begin{equation}\label{dist}
   \textrm{dist}(x, \partial \Omega) \le (\frac{\varepsilon}{C_{i}})^{1/ \beta_{i}}  \quad \forall i=1, \cdots ,m.
      \end{equation}
    Let $y\in \partial \Omega$ be a point such that has minimum distance to $x$. Then  by assumption on the boundary values, there is $k$ such that $u_{k}(y)=0$ and
    \[
    \frac{ |u^{\varepsilon}_{k}(x)- u^{\varepsilon}_{k}(y) |   }{|x-y |^{\beta_{k}}} \le C_k.
    \]
    The previous inequality  and  (\ref{dist})  imply that
    \begin{equation}\label{bound2}
    u^{\varepsilon}_{k}(x) \le \varepsilon.
    \end{equation}
Combining (\ref{bound1}) and (\ref{bound2}) yields that Laplace of $u_i$ is bounded. 



\begin{rem}
The uniform bound of  normal derivative of  $ u_{i}^{\varepsilon} $ yields estimates for limiting problem as follows.
Integrate from
 \[
 \Delta  u_{i}^{\varepsilon}=  \frac{A_i(x) }{\varepsilon}   \prod\limits_{i=1}^{m}  u_{i}^{\varepsilon}  (x)
 \]
 to obtain
\[
\int_{\partial\Omega}  \frac{\partial u_{i}^{\varepsilon}}{\partial n}\, ds =  \int_{\Omega} \frac{A_i(x) }{\varepsilon}   \prod\limits_{i=1}^{m}  u_{i}^{\varepsilon}  (x)\,  dx.
\]
From here we get
 \[
  \int_{\Omega} \frac{A_i(x) }{\varepsilon}   \prod\limits_{i=1}^{m}  u_{i}^{\varepsilon}  (x)\,  dx \le C
  \]
  this shows
  \[
  \int_{\Omega} A_i(x)   \prod\limits_{i=1}^{m}  u_{i}^{\varepsilon}  (x)\,  dx    \rightarrow 0 \quad \text{ as $ \varepsilon$ tends to zero. }
  \]
\end{rem}
\begin{defn}
Consider  the non empty  sets    $\Omega_{i}:={\{ x \in \Omega : u_{i}(x) =0 }\}$. Then the  free boundaries (interfaces)  are define as
\[
\Gamma_{i,j} = \partial\Omega_{i} \cap  \partial \Omega_{j}\cap \Omega.
\]
\end{defn}
In the next Lemma we give the free boundary condition for the case $A_{i} =1$.
 \begin{lem}\label{F3}
 The following conditions holds on  the free boundary  $\Gamma_{i,j}.$
\begin{enumerate}
        \item   $
 \frac{\partial u_i}{\partial n}|_{\Omega_j} =-  \frac{\partial u_j}{\partial n}|_{\Omega_i},\\
 $
 \item
 $
  \frac{\partial u_k}{\partial n}|_{\Omega_j} -  \frac{\partial u_k}{\partial n}|_{\Omega_i}=\frac{\partial u_i}{\partial n}|_{\Omega_j}  \quad    k\neq i,j.
  $
  \end{enumerate}
 \end{lem}
 \begin{proof}
 Let $x_0$ be a free boundary point in $\Gamma_{i,j}.$
 Note that
 \[
 \Delta (u_k -u_j) =0  \quad \text{in } B_{r}(x_0) \setminus \Gamma_{i,j}.
 \]
 In the sense of distribution we have
 \[
 \Delta (u_k -u_j) = \frac{\partial (u_k - u_j)}{\partial n} H^{n-1}|_{\Gamma_{i,j}} \quad \text{in  } \, B_r,
 \]
 Splitting  $B_r= (B_{r} \cap   \Omega_i) \cup (B_{r} \cap   \Omega_j) $ and considering the fact that in $\Omega_j $ we have  $u_j = 0 $ the second relation is proved.

 \end{proof}

\begin{rem}:
 In  \cite{W}   the  uniqueness of the limiting solution   of system (\ref{f20})  for  arbitrary number of components, is shown.  Consider  the metric space  $ \sum $ defined by
\[
\sum = {\{ (u_1, u_2 , \cdots, u_m) \in \mathbb{R}^m :  u_{i}\ge 0,\, u_{i}\cdot u_{j}=0 \quad \text{for} \quad i\neq  j\}}.
    \]
In  \cite{W}  (see Theorem 1.6 ) it is   shown  that    the limiting   solution    $(u_1, \cdots ,u_m)$   of   (\ref{f20})   is a harmonic map into the space  $ \sum.$      By definition the harmonic map is the critical point  of the following energy  functional
	\[
	\int_{\Omega}  \sum_{i=1}^{m}  \frac{1}{2}| \nabla u_{i}|^{2} dx,
	\]
	among all nonnegative  segregated states   $u_i \cdot  u_j = 0,$  a.e.  with the same boundary conditions.

 Also in \cite{AB} an alternative proof of uniqueness for limiting case  for system(\ref{f20}) is given  which   is more direct and based on properties of limiting solutions. Although
some  properties of limiting solution for systems (\ref{s1}) and  (\ref{f20})  are similar, the proof of  uniqueness  for system   (\ref{s1}) in the case  $\varepsilon$ tends to zero remains challenging  problem.
\end{rem}
Define the energy associated to  $m$ densities defined by
\begin{equation*}
E(U)=\int_{\Omega} \sum _{ i}|\nabla u_{i}(x)|^{2} dx.
\end{equation*}
Now consider the  following problem
\[
\min  E(U),
\]
over the closed but non-convex  set
\begin{equation*}
S= \left\{ (u_1, \cdots u_m): u _{i} \in H^{1}(\Omega),\, u_{i} \geq 0, \,  \prod\limits_{i=1}^{m}  u_{i} (x)=0, \, u_{i}|_{\partial\Omega}=\phi_{i} \right\}.
\end{equation*}

  Existence of a minimizer is direct. The following variation
\[
v_{i} = (1 + \varepsilon  \varphi_{i}) u_i, \, i=1,\cdots,m,
\]
with $ \varphi_{i} \in C^{\infty}_{c}(\Omega)$ yields the following
\[
u_{i}\ge 0,  \quad  u_{i}\cdot \Delta u_{i} =  0, \quad  \Pi_{j=1}^{m} u_j=0.
\]
This implies that  each  $u_i$ is harmonic in its support which dose not hold  for our limiting solution.  In fact,  for system (\ref{s1}) in Theorem 3.2 we show that
\[
 \prod\limits_{i=1}^{m}  u_{i} =0 \quad \textrm{and} \, \Delta u_{i}  \, \textrm{is bounded}.
  \]
   Figure 1 also shows $u_1$ is not smooth in its support    $ \Delta u_{i}  $  are Dirac measures on interfaces.

\section{Explicit solutions in the limiting case}\label{explicit}

In  this section  we give an explicit solution and the rate of convergence for  the limiting solution of the following system
\begin{equation}\label{s113}
\left \{
\begin{array}{lll}
\Delta  u_{i}^{\varepsilon}=   \frac{ A_{i} (x) }{\varepsilon}   \prod\limits_{j=1}^{m}  (u_{j}^{\varepsilon})^{\alpha_j}   (x) & \text{ in  } \Omega,\\
u_{i} =\phi_{i}  &   \text{ on   } \partial \Omega,
 \end{array}
\right.
\end{equation}
 for the cases that  $A_{i}(x)$ are the same or $A_i=C_i,$ constants.

\subsection{Construction of Solutions}

It is easy to check that for every $ \varepsilon$
\[
 \Delta( u_{1}^{\varepsilon}-  u_{i+1}^{\varepsilon}) =0,  \quad i=1,\cdots ,m-1
 \]
 which   remains  true as  $\varepsilon$ tends to zero. First of all define
\begin{equation}
	w_i = u_1-u_{i+1}, \qquad i=1,\ldots,m-1,
\end{equation}
then $w_i$ is the harmonic extension of the Dirichlet value
 $\phi_1 - \phi_{i+1}.$  This means that $w_i$  for $i=1, \cdots, m-1$ is the solution of
 \begin{equation}\label{sunshine3}
\left \{
\begin{array}{lll}
\Delta w_{i}= 0   & \text{ in  } \Omega,\\
w_{i}  =\phi_{1}-  \phi_{i+1}    &   \text{ on   } \partial \Omega.
 \end{array}
\right.
\end{equation}

Note that the nonnegativity of the $u_i$ is equivalent to $u_1 \geq w_i$. Thus, an obvious candidate solution
is given by
\begin{equation}\label{m1}
	u_1(x) = \max \,  \left( \underset{i=1,\ldots,m-1}{\max } w_i(x), 0 \right)\\
\end{equation}
and
\begin{equation}\label{m2}
	u_i = u_1 - w_i, \qquad  i=2,\ldots,m.
\end{equation}
Obviously, by this construction we have $u_i \geq 0$ and moreover
 \[
 u_1(x) u_2(x) \ldots u_m(x) = 0, \quad \text{  for\,  all  } \quad x \in \Omega.
 \]
 To see the latter, let $x$ be fixed and $j$ such that $w_j(x) \geq w_i(x)$ for all $i$. Then
 \[
 u_j(x) = u_1(x) - w_j(x) = 0.
 \]
 We finally need to verify $\Delta u_i \geq 0$. For $u_{1}$  this follows from the fact that maximum of harmonic function is subharmonic then for the rest of $u_{i}$ it follows from (\ref{m1}) and  (\ref{m2}).

\begin{rem}
 Let  $v_{i}$  be defined as  below :
\[
v_i=u_{2}- u_{i},
\]
then set
\[
 u_{2}= \max{\{\max_{i=1,\ldots,m-1} v_i(x), 0}\}.
 \]
 From this we can recover other components by
 \[
 u_{i}= v_i-u_{2}, \quad i=1,3,\cdots m.
\]
One can check this choice  gives the same solutions as in (\ref{m1})
and (\ref{m2}), for the case $m=3$ is straightforward.
\end{rem}

\subsection{Convergence Rate}

We now turn our attention to a rate of convergence of the solutions as $\varepsilon \rightarrow 0$. Note that
\begin{equation}
	w_i = u_1^\varepsilon-u_i^\varepsilon, \qquad i=1,\ldots,m
\end{equation}
is harmonic with Dirichlet data $\phi_1 - \phi_i$, hence coincides with the one in the previous section, in particular independent of $\varepsilon$.

We thus have
$$
\Delta  u_{1}^{\varepsilon}=   \frac{ 1 }{\varepsilon}   \prod\limits_{i=1}^{m}  u_{i}^{\varepsilon}
=  \frac{ 1 }{\varepsilon}   \prod\limits_{i=1}^{m}  (u_{1}^{\varepsilon}-w_i)
$$
Now we have $0 \leq u_{1}^{\varepsilon}-w_i$ and $u_{1}^{\varepsilon}-w_i \geq u_1^\varepsilon - u_1$, hence
$$\varepsilon \Delta  u_{1}^{\varepsilon} \geq  |u_1^\varepsilon - u_1|^m, $$
respectively
$$ \varepsilon \int_\Omega |\nabla (u_1^\varepsilon - u_1)|^2 ~dx + \int_\Omega |u_1^\varepsilon - u_1|^{m+1} ~dx \leq
 - \varepsilon \int_\Omega \nabla (u_1^\varepsilon - u_1)\cdot \nabla u_1 ~dx .$$
Applying Young's inequality on the right-hand side we deduce
$$ \Vert u_1^\varepsilon - u_1 \Vert_{L^{m+1}(\Omega)} \leq C \varepsilon^{1/(m+1)}. $$


\section{Numerical Study of the Limiting Problem}
 This section provides some examples of numerical approximations to the limiting problem of the  following
 system.
\begin{equation}\label{square1}
\left \{
\begin{array}{lll}
\Delta u_{i}= \frac{1}{\varepsilon}   \prod\limits_{j }^{m}  \,  u_{j}  & \text{ in }\ \Omega, \\
   u_{i}  =\phi_{i}   &  \text{ on } \  \partial\Omega.
\end{array}
\right.
\end{equation}
In our examples we implemented     directly mimicking the fixed point technique in the existence proof of Theorem \ref{sun0} with value of $\varepsilon$  and the method in Section 4 which  demonstrate those give basically the same as epsilon goes to zero

\begin{exam}\label{E1}

Let $\Omega =B_{1}, m=3.$   The boundary values $ \phi_{i}$ for  $ i=1,2,3$ are defined by
\begin{equation*}
\phi_{1}(1,\Theta)=
\left \{
\begin{array}{ll}
|\sin(\frac{3}{2}\Theta)|  &  0 \leq\Theta \leq \frac{4\pi}{3},\\
0 & \text{ elsewhere,}
\end{array}
\right.
\hspace{0.1in}
\phi_{2}(1,\Theta)=
\left \{
\begin{array}{ll}
|\sin(\frac{3}{2}\Theta)|  &   \frac{2\pi}{3} \leq\Theta \leq 2\pi,\\
0 & \text{elsewhere.}
\end{array}
\right.
\end{equation*}
 \begin{equation*}
\phi_{3}(1,\Theta)=
\left \{
\begin{array}{lr}
|\sin(\frac{3}{2}\Theta)|  &   \frac{4\pi}{3} \leq\Theta \leq 2\pi + \frac{2\pi}{3},\\
0 &  \text{elsewhere.}
\end{array}
\right.
\end{equation*}
Here the boundary conditions satisfy
\[
\phi_{1}\cdot \phi_{2} \cdot \phi_{3}= 0.
\]
 The surface of $u_1$ is depicted in Figure 1. Also one can check the jump in gradient  of  $u_1$ along $ \Gamma_{2,3} $
 which has shown in part  2 of  Lemma (\ref{F3}).  In Figure 2,

\begin{figure}[h!]
{\includegraphics[width=.7 \columnwidth]{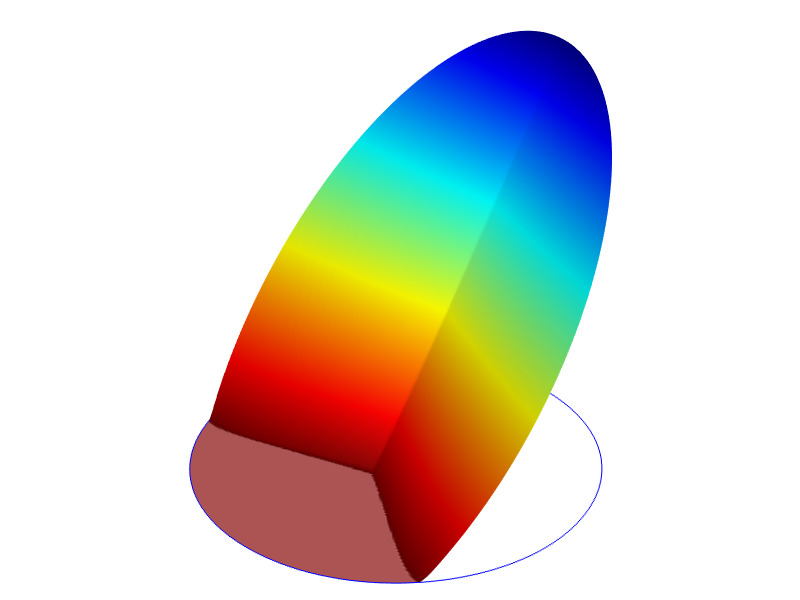}}
  \caption{ surface of $u_1$ }
  \label{fig:E}
\end{figure}

\begin{figure}[h!]
{\includegraphics[width=.6 \columnwidth]{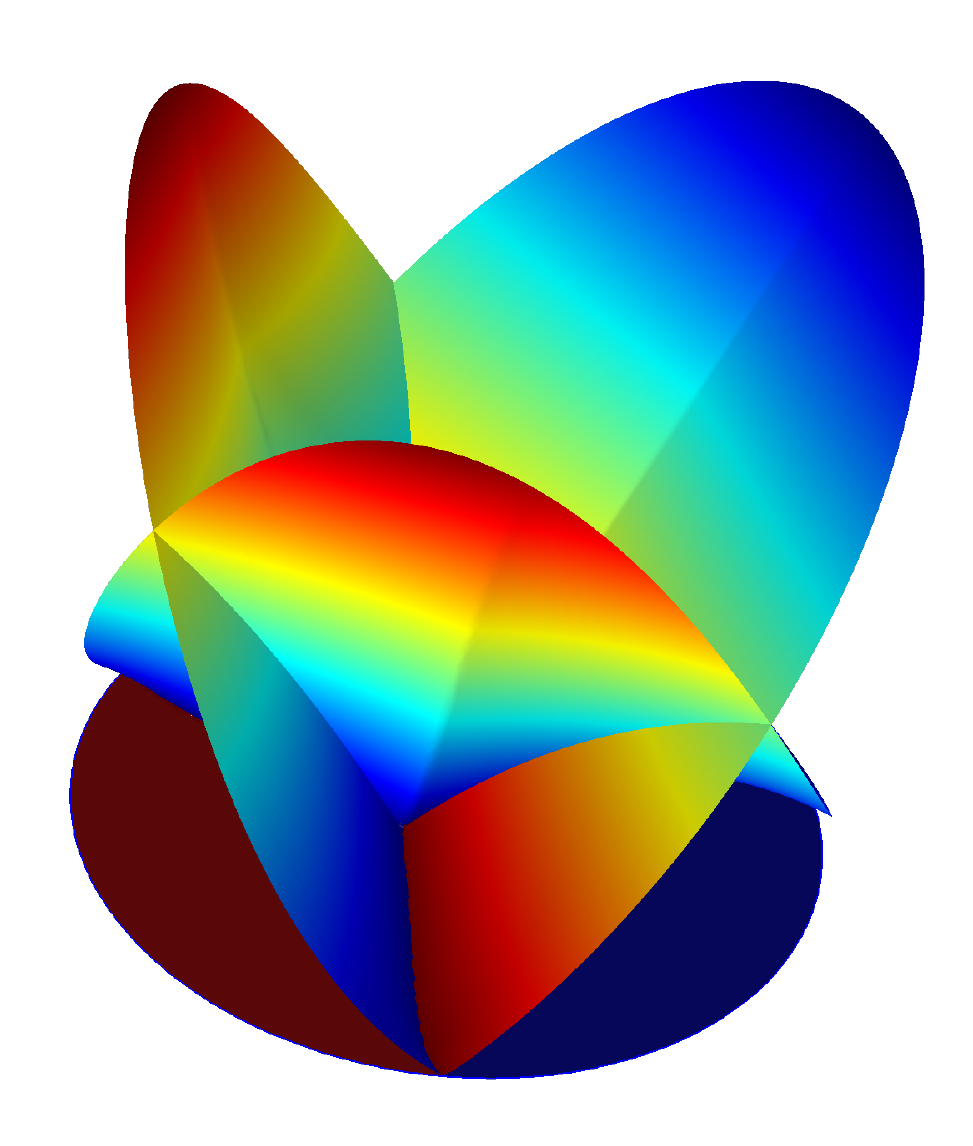}}
  \caption{surface of $u_1 + u_2 + u_3$.}
  \label{fig:E}
\end{figure}

\end{exam}

\begin{exam}\label{farid20}
Let $\Omega =[-1,1]\times [-1, 1]$ and $m=4.$
The boundary values $\phi_{i},$ (i=1,2,3,4) are given as follows:
\begin{equation*}
\phi_{1}=
\left \{
\begin{array}{lr}
1-x^{2} &   x\in[-1, 1] \  {\&} \  y=1 ,\\
0 & \ \ \text {elsewhere.}
\end{array}
\right.
\hspace{0.1in}
\phi_{2}=
\left \{
\begin{array}{lr}
2(1-y^{2}) & y \in [-1, 1] \  {\&} \  x=1 ,\\
0 & \ \ \text {elsewhere.}
\end{array}
\right.
\end{equation*}

\begin{equation*}
\phi_{3}=
\left \{
\begin{array}{lr}
3(1-x^{2}) &  x\in [-1 ,1] \  {\&} \  y=-1 ,\\
0 & \ \ \text {elsewhere.}
\end{array}
\right.
\hspace{0.1in}
\phi_{4}=
\left \{
\begin{array}{lr}
4(1-y^{2}) &   y \in [-1 , 1] \  {\&} \  x=-1 ,\\
0 & \ \ \text {elsewhere.}
\end{array}
\right.
\end{equation*}
We implemented the iterative scheme given by Lemma 2.2  with $ \varepsilon=10^{-8}$ and method given by (\ref{m1}) and (\ref{m2}). The obtained solutions are same and the surface of $u_1$ is given in  (\ref{fig:EL}).
\begin{figure}[h!]
{\includegraphics[width=.8 \columnwidth]{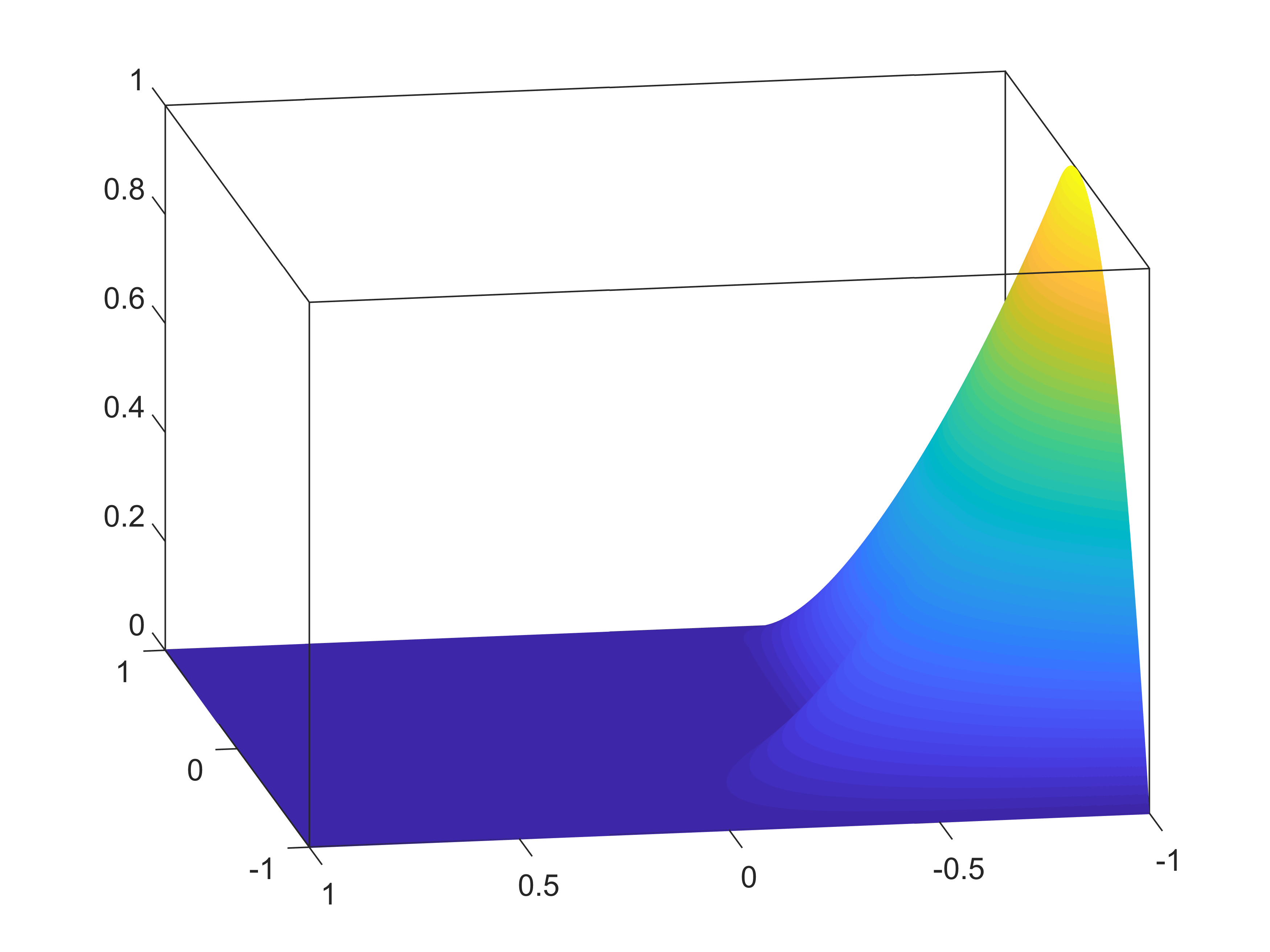}}
  \caption{surface of $u_1$.}
  \label{fig:EL}
\end{figure}
The interfaces are shown in Figure \ref{fig5}.


\begin{figure}[!hbt]
\begin{picture}(280,280)(0,0)
\put(0,0){\includegraphics[height=280pt]{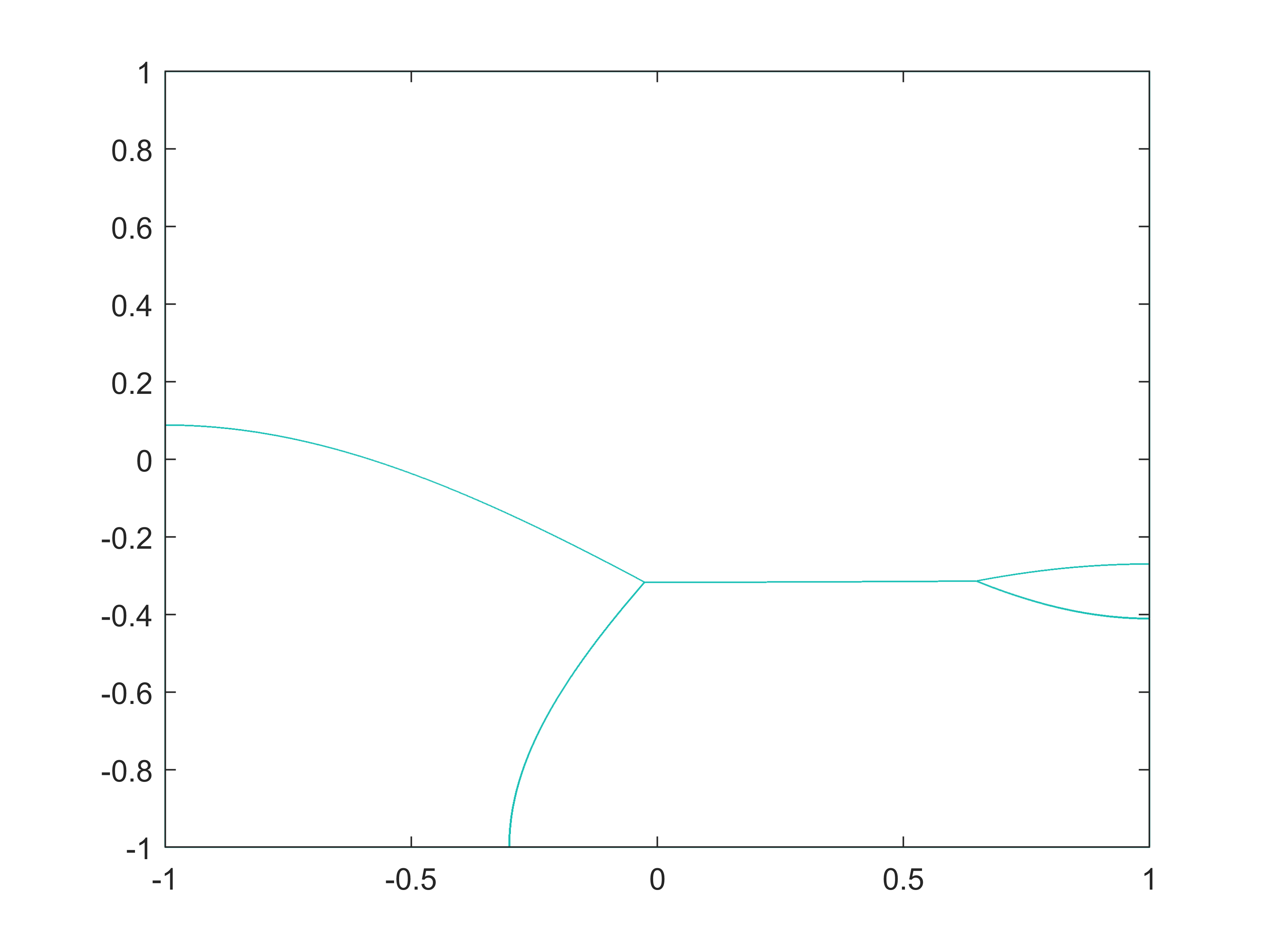}}
  \put(160,180){$u_{1}=0 $}
\put(295,105){$u_{4}=0$ }
\put(200,70){$u_{3}=0$ }
\put(100,70){$u_{2}=0$ }
\end{picture}
\caption{Free boundary and supports of the components.}
\label{fig:D}
\end{figure}

In Figure (\ref{fig6} we draw the Laplace of $u_1$ on the interfaces.  We know Laplace of $u_1$ is Dirac along interfaces so we  scaled   $\Delta u_1$  by multiplying by mesh size.

\begin{figure}[h!]
{\includegraphics[width= \columnwidth]{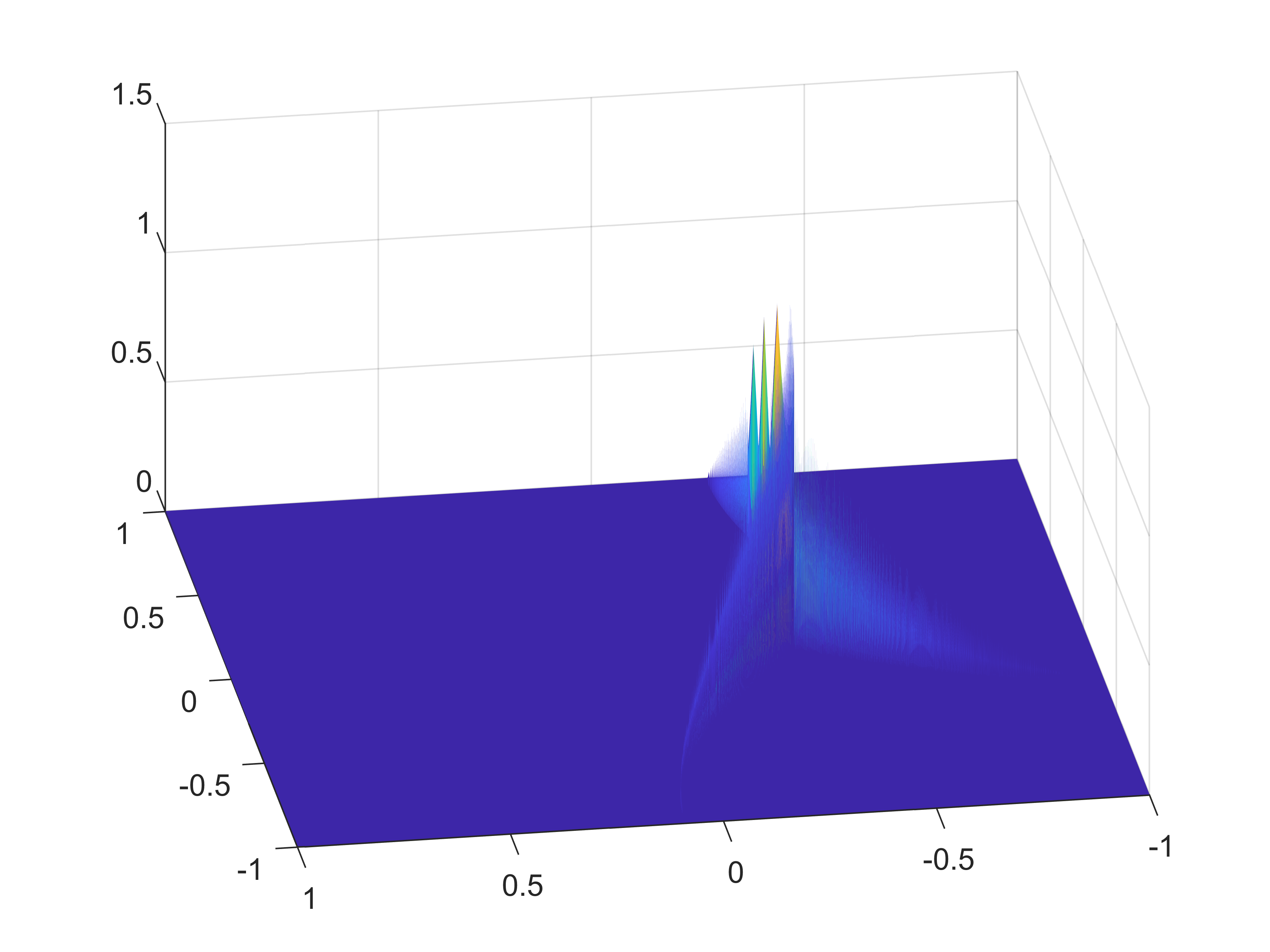}}
  \caption{ Laplace of $u_1$ as measure(scaled) on the interfaces. The mesh size is $\triangle x=\triangle y=10^{-3}.$}
  \label{fig6}
\end{figure}
\end{exam}

\begin{exam}\label{farid22}
Next, we change boundary  values as below.
\begin{equation*}
\phi_{1}=
\left \{
\begin{array}{lll}
1-x^{2} &   x\in[-1, 1] \  {\&} \  y=1 ,\\
1-y^{2} &   y\in[-1, 1] \  {\&} \  x=1 ,\\
0 & \ \ \text {elsewhere.}
\end{array}
\right.
\hspace{0.1in}
\phi_{2}=
\left \{
\begin{array}{lll}
2(1-y^{2}) & y \in [-1, 1] \  {\&} \  x=-1 ,\\
2(1-x^{2}) & x \in [-1, 1] \  {\&} \  y=1 ,\\
0 & \ \ \text {elsewhere.}
\end{array}
\right.
\end{equation*}

\begin{equation*}
\phi_{3}=
\left \{
\begin{array}{lll}
3(1-x^{2}) &  x\in [-1 ,1] \  {\&} \  y=1 ,\\
3(1-y^{2}) &  y\in [-1 ,1] \  {\&} \  x=1 ,\\
0 & \ \ \text {elsewhere.}
\end{array}
\right.
\hspace{0.1in}
\phi_{4}=
\left \{
\begin{array}{lll}
4(1-x^{2}) &   -1\le x \le 1 \  {\&} \  y=-1 ,\\
4(1-y^{2}) &   y \in [-1 , 1] \  {\&} \  x=-1 ,\\
0 & \ \ \text {elsewhere.}
\end{array}
\right.
\end{equation*}

  The following picture shows the interfaces

\begin{figure}[h!]
{\includegraphics[width=.8 \columnwidth]{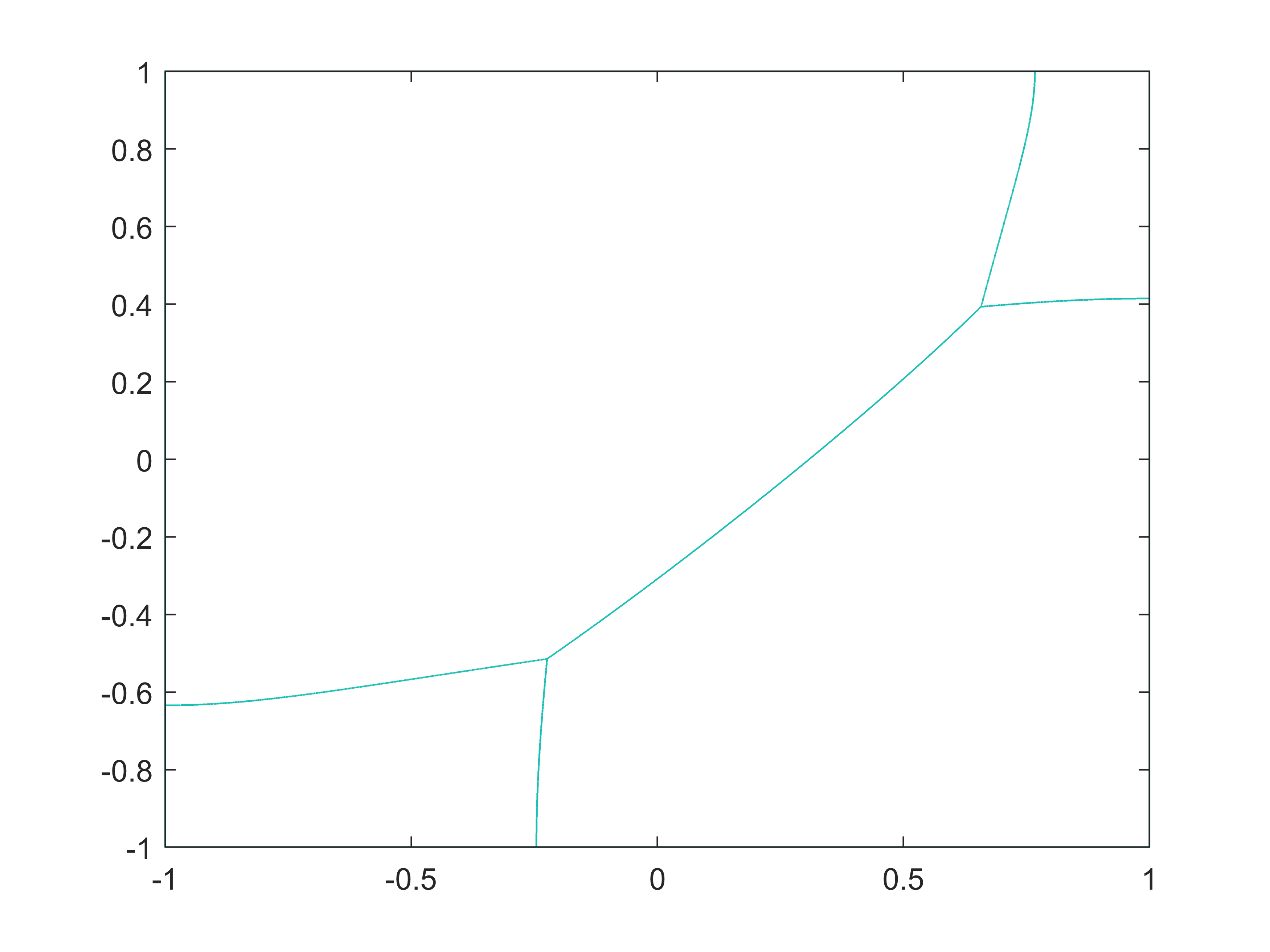}}
  \caption{Free boundaries.}
  \label{fig:EL}
\end{figure}
%

\end{exam}

\renewcommand{\refname}{REFERENCES }

\end{document}